\documentclass[10pt,fleqn]{article}

\usepackage{amsthm,amsfonts,amssymb,indentfirst,latexsym,bm,amsmath}

\usepackage{amsmath}

\usepackage{amsmath}
\usepackage{amsfonts,amsthm,amssymb}
\usepackage{amsfonts}
\usepackage{graphics}
\usepackage{mathrsfs}

\begin{document}
\title{ Concentration properties of semi-vertex transitive graphs and random bi-coset graphs\thanks{
The project supported partially by CNNSF (No.10971106).}}
\author{Xingchao Deng$^a$  and Kainan Xiang$^b$\\
 \footnotesize $^a$Center for Combinatorics, LPMC-TJKLC\\
 \footnotesize Nankai University, Tianjin 300071, P. R. China\\
 \footnotesize  E-mail: dengyuqiu1980@126.com\\
   \footnotesize $^b$School of Mathematical Sciences, LPMC, Nankai University\\
       \footnotesize Tianjin City, 300071, P. R. China\\
       \footnotesize E-mail: kainanxiang@gmail.com}
\date{}
\maketitle
\begin{abstract}
It is well-known that concentrators are sparse graphs of high
connectivity, which play a key role in the construction of switching
networks; and any semi-vertex transitive graph is isomorphic to a
bi-coset graph. In this paper, we prove that random bi-coset graphs
are almost always concentrators, and construct some examples of
semi-vertex transitive concentrators.

\end{abstract}
\section{Introduction}
Many problems on information transmission and complexity have shown
the importance of constructing graphs that are highly connected yet
sparse. On one hand, concentrators are a key building block for the
construction of a class of graphs called superconcentrators which
are useful in the study of algorithmic complexity. On the other
hand, concentrators are the basis for another class of graphs called
generalized connectors($\cite{20}$). For the importance of
concentrator-like bipartite graphs in constructing low complexity
error-correcting codes, refer to Tanner$\cite{19}.$

  Recall a bipartite graph is an $(n,\theta,k,\alpha, c)$ bounded strong concentrator
(bsc) if it is a bipartite graph with $n$ inputs, $\theta n$ outputs
and at most $kn$ edges such that $\Gamma_X \geq c\vert X\vert$ for
any set $X$ of inputs with $\vert X\vert \leq \alpha n .$ Here
$\Gamma_X$ is the set of outputs connected to $X$ and
$\vert\cdot\vert$ is the cardinality of a set.

   An $(n,k)$-superconcentrator is a directed graph with $n$ inputs and $n$ outputs,
and at most $kn$ edges satisfying that for any $1\leq r\leq n$ and
any two sets of $r$ inputs and $r$ outputs, there are $r$ vertex
disjoint paths connecting the two sets. A family of linear
superconcentrators of density $k$ is a set of $(n_i,k+o(1))$
superconcentrators, with $n_i\rightarrow \infty$, as $i\rightarrow
\infty$, which is most useful in theoretical computer science. Note
an $(n,k,c)-expander$ is a bipartite graph with $n$ inputs, $n$
outputs and at most $kn$ edges such that for any subset $A$ of
inputs,
 $$\vert N(A)\vert \geq \left(1+c\left(1-\frac{\vert
    A\vert}{n}\right)\right)|A|,$$
where $N(A)$ is the set of all neighbors of $A.$

A graph $G=(V,E)$ on $n$ vertices with maximal degree $d$ is called
an $(n,d,c)-magnifier$ if $\vert N(X)-X\vert\geq c\vert X\vert $ for
any vertex set $\vert X\vert\leq \frac{n}{2}.$ For a graph $G=(V,E)$
with $V=\{v_1,v_2,\cdots ,v_n\},$ its (extended) double cover is the
bipartite graph $H$ on input set $X=\{x_1,\cdots,x_n\}$ and output
set $Y=\{y_1,\cdots,y_n\}$ such that $x_i \in X $ and $y_j \in Y$
are adjacent if and only if $i=j$ or $v_iv_j\in E.$\\

\noindent{\bf Lemma $1.1^{\cite{1}}.$} Let graph $G=(V,E)$ be an
$(n,d,c)-magnifier$ and $H$ its extended double cover. Then $H$
is an $(n,d+1, c)-expander.$\\

By this lemma, to construct a family of linear expanders, it
suffices to construct a family of linear magnifiers, i.e., to
construct, for some fixed $d$ and $c>0$, a family of
$(n_i,d,c)-magnifiers$ with $n_i\rightarrow \infty$.

It is well-known that random $d-regular$ graphs are magnifiers for
$d\geq 5.$ But an explicit construction is needed for application.
However, such construction is much more difficult. Margulis gave the
first explicit family of linear magnifiers of density 5 and proved
it has expansion $c$ for some $c>0$ by several deep results from the
theory of group representations. But he didn't bound $c$ strictly
away from 0. Then Gaber and Galil modified Margulis' construction to
obtain a family of linear magnifiers with density 7 and expansion
$\frac{2-\sqrt{3}}{2};$ and used this family to construct explicitly
a family of linear superconcentrators of density $271.8.$
Sch\"{o}ing $\cite{15}$ constructed the smaller superconcentrators
of density 28, which is the best density.

Note the Cayley diagram of a group $G$ with respect to a multiset
$S$ of elements of $G$ is a directed multigraph on $G$ whose
multiple edges are all ordered pairs $(x,y)\in G^2$ with $y=sx$ for
some $s\in S.$ Then Cayley graph $X(G,S)$ of $G$ with respect to $S$
is the Cayley diagram with ignoring orientation but remaining
multiple edges.

Recall many  constructions of magnifier graphs are Cayley graphs,
and many families of finite simple groups are magnifier families.
Kassabov $\cite{9}$ constructed explicitly generating sets $S_n$ and
$\widetilde{S}_n$ of the alternating and the symmetric groups
respectively to obtain two families $\left\{X(Alt(n),
S_n)\right\}_{n\geq 1}$ and $\left\{X\left(Sym(n),
\widetilde{S}_n\right)\right\}_{n\geq 1} $ of bounded degree
magnifiers of Cayley graphs. Kassabov et al $\cite{10}$  proved that
there exist $k\in \mathbb{N}$ and $0 < \varepsilon <\infty$ such
that any non-abelian finite simple group $G$, which is not a Suzuki
group, has a set of $k$ generators for which the Cayley graph $X(G,
S)$ is an $\varepsilon-expander$ (here we call
$\varepsilon-magnifier$). Notice we will use the result in [9] to
construct concentrators in Section 6.

For any graph $H,$ let $\mu_1(H)$ denote the second largest
eigenvalue in absolute value of its adjacency matrix $A_H.$ When $H$
is $d$-regular, the normalized adjacency matrix
$A^*_H=\frac{1}{d}A_H$ of $H$ is doubly stochastic and $\mu_1^*(H)=
\frac{1}{d}\mu_1(H);$ where $\mu_1^*(H)$ is the second largest
eigenvalue in absolute value of $A^*_H.$\\

\noindent{\bf Theorem (Alon-Roichman) $1.2^{\cite{4}}.$} For any $\varepsilon
>0,$ there is a $c(\varepsilon)>0$ depending only on $\varepsilon$
such that the following holds. Let $G$ be a group of order $n$ and
let $S$ be a set of $c(\varepsilon ) \log n$ elements of $G$ chosen
uniformly and independently at random. Then
$$\mu_1^*(X(G,S))\ \mbox{is\
  concentrated\ around\ its\ mean\ and}\
  E(\mu_1^*(X(G,S)))<\varepsilon .$$\\

By Theorem 1.2, we have\\

\noindent{\bf Corollary $1.3^{\cite{6}}.$}  For any $\varepsilon
 >0,$ there is a $c^*(\varepsilon) > 0$ depending only on
$\varepsilon$ satisfying that for any finite group $G$ with $n$
elements, the Cayley graph $X(G, S)$ is an $(n, 2|S|,
\varepsilon)-expander\ (magnifier)$ with high probability as
$n\rightarrow \infty$, where $S$ is a multiset of
$c^*(\varepsilon) \log n $ random elements of $G$.\\

Christofides  and  Markstr\"{o}m$\cite{6}$ generalized Alon-Roichman
theorem to random coset graphs. For any $p\in (0,1),$ let $H_p$ be
the $weighted\ entropy\ function$ defined by
$$H_p(x) = x \log \frac{x}{p}+ (1-x) \log \frac{1-x}{1-p},\ x\in [0,1].$$
Then by [6], the following theorem holds.\\

\noindent{\bf Theorem $1.4^{\cite{6}}.$} Let $S$ be a multiset of
$k$ elements of a finite group $G$ chosen independently and
uniformly at random, and $D$ the sum of the dimensions of the
irreducible representations of the group $G.$ Then for any $0
<\varepsilon < 1,$
$$P_r(\mu^*(X(G, S)) > \varepsilon)\leq 2D\exp \left\{-k H_ \frac{1}{2}\left(\frac{1+\varepsilon}{2}\right) \right\}.$$\\

For any subgroup $H$ of $G,$ denote by $D(G,H)$ the sum of the
dimensions of the irreducible representations of $G$ which do not
contain the trivial representation of $H$ when decomposed into
irreducible representations of $H.$

Let $S$ be a multiset of $G$  chosen independently and uniformly at
random. The random coset graph $X(G,H;S)$ of $G$ with respect to $H$
and $S$ is defined as follows: its vertices are all right cosets of
$H$ in $G,$ and there is an edge between $Hg_1$
and $Hg_2$ if and only if $g_2g_1^{-1}\in HSH.$\\

\noindent{\bf Theorem $1.5^{\cite{6}}.$} For the random coset graph
$X=X(G,H;S)$ and any $0 <\varepsilon < 1,$
 $$P_r(\mu^{*}(X)>\varepsilon) \leq2D(G,H)
    \exp\left\{ -|S| H_ \frac{1}{2}\left(\frac{1+\varepsilon}{2}\right)\right\}.$$\\

Note the extended double cover is a superconcentrator. Moreover superconcentrators can be
constructed by bounded concentrators through recursive construction of$\cite{21}$. The Cayley graph is vertex
transitive. Tanner[18] constructed several explicit concentrators by
the generalized polygons and his technique for deciding the
concentration properties of a graph by analysing its eigenvalues can
be stated as follows:

Assume $I$ and $O$ are two disjointed sets of sizes $n$ and $m$
respectively. Let $H$ be a bipartite graph with $I$ as input vertex
set and $O$ as output vertex set such that edges connect input
vertices to output ones, and the degree of each input vertex is $k$
and that of each output vertex is $r=n\frac{k}{m}.$ Write
$A=[a_{ij}]$ for the incidence matrix of $H:$ $a_{ij}=1$ if the
$ith$ input vertex $a_i$ is connected to the $jth$ output vertex
$b_j$ and $=0$ otherwise. Note $AA^T$ is diagonalizable and has real
nonnegative eigenvalues due to it is symmetric and nonnegative
definite. If $\lambda_1\geq \lambda_2\geq \cdots\geq \lambda_n$ are
the ordered eigenvalues of $AA^T,$ then the following holds.\\

\noindent{\bf Theorem $1.6^{\cite{18}}.$} If $\lambda_1>\lambda_2,$
then for any $0<\alpha\leq \frac{m}{n},$ $H$ is an
$\left(n,\frac{m}{n},l,\alpha,c(\alpha)\right)~bsc$ with
$$ c(\alpha)\geq\frac{l^2}{[\alpha(kr-\lambda_2)+\lambda_2]}.$$\\

\noindent{\bf Definition $1.1.$} For a bipartite graph $H$ with the
vertex bipartition $U(H)$ and $W(H),$ $H$ is called semi-vertex
transitive if $Aut(H)$ is transitive on $U(H)$ and $W(H)$
respectively.\\

\noindent{\bf Definition $1.2.$} Let $G$ be a group with $L$ and $N$
as its two subgroups. For a set $S$ of some bi-cosets $NgL,$ define
the bi-coset bipartite graph $X=C(G,L,N;S)$ of $G$ with respect to
$L,N,S$ as follows: its vertex set $V(X)$ is $[G:L]\cup [G:N],$ and
its edge set $E(X)$ is
 $$\{\{Lg,Nsg\}:g\in G,s\in S\}.$$
Particularly, $BC(G,S):=C(G,\{1\},\{1\};S)$ is called the bi-Cayley
graph of $G$ with respect to $S.$ In addition, if $S$ is a multiset
of $G$ chosen independently and uniformly at random, then we
call $X=C(G,L,N;S)$ a random bi-coset graph of $G$ with respect to $L,\ N$ and $S.$\\

The following Proposition 1.7 on bi-coset bipartite graphs is well-known:\\

\noindent{\bf Proposition $1.7.$} $\mathbf{(i)}$ Let $X=C(G,L,N;S)$
be the bi-coset graph of $G$ with respect to $L,\ N,\ \mbox{and}\ S.
$ Then the following hold.
\begin{eqnarray*}
&&\bullet\ \mbox{The\ degrees\ of}\ Lg\ \mbox{and}\ Ng\
   \mbox{are}\ |S:N|\ \mbox{and}\ |S:L|\ \mbox{respectively\ for\ any}\ g\in G,\\
&&\ \ \ \mbox{and}\ X\ \mbox{is\ regular\ if\ and\ only\ if}\
     |L|=|N|.\\
&&\bullet\ G\subseteq Aut(X)\ (\mbox{acting\ on}\ [G:L]\
    \mbox{and}\ [G:N]\ \mbox{by\ right\ multiplication}),\\
&&\ \ \ X \ \mbox{is\ semi-vertex\ transitive.}\\
&&\bullet\  X\ \mbox{is\ connected\ if\ and\ only\ if}\
              G=\langle S^{-1}S\rangle.
\end{eqnarray*}
$\mathbf{(ii)}$ Every semi-vertex transitive graph is isomorphic to
some bi-coset graph.\\

For any bi-coset graph $X=C(G,L,N;S)$ of $G$ with respect to $L,N$
and $S,$ let $A=[a_{ij}]$ be its incidence matrix and
$M_L=\frac{1}{2|S|^2}AA^T.$ Then one of our main results is stated
as follows.\\

\noindent{\bf Theorem $1.8.$} Let $X=C(G,L,N;S)$ be a random
bi-coset graph of $G$ with respect to $L,N$ and $S\ (|S|=k).$ Then
for any $0<\varepsilon<1,$
$$P_r\left(\mu^* \left(M_L\right)>\varepsilon+\frac{1-\varepsilon}{2k}\right)\leq
  2D\left(G,L\right)exp\left\{-(k^2-k) H_{\frac{1}{2}}\left(\frac{1}{2}+\varepsilon  \right)\right\}.$$\\

By Theorem 1.8 and Proposition 1.7(ii), we have\\

\noindent{\bf Corollary $1.9.$} The set of all semi-vertex
transitive graphs isomorphic to some bi-coset graph
$C\left(G,L,N;\widehat{S}\right)$ with
$\left\vert\widehat{S}\right\vert =k$ is just the sample space for
the random bi-coset graph $X=C(G,L,N;S)$ with $\vert S\vert =k;$ and
hence in this sense, almost always semi-vertex
transitive graphs are concentrators.\\

Our other main results are constructions of semi-vertex transitive
concentrators which are presented in Sections 4-6, which may be simple but we can not find them in the former papers. Moreover we give the theory explanation of the golay codes constructed by Mathieu groups are good.

\section{Preliminaries from representation theory}
A $representation~\rho$ of a finite group $G$ is a homomorphism
$\rho :G\rightarrow \bigcup(\mathbb{H}),$ where $\mathbb{H}$ is a
finite dimensional Hilbert space and $\bigcup(\mathbb{H})$ is the
group of unitary operators on $\mathbb{H}.$ The dimension $d_\rho$
of $\rho$ is the dimension of $\mathbb{H}.$ Fix a basis for
$\mathbb{H},$ then each $\rho(g)$ is associated with a unique
unitary matrix $[\rho(g)]$ satisfying
$[\rho(gh)]=[\rho(g)][\rho(h)]$ for any $g,h \in G.$

For a fixed representation $\rho :G\rightarrow \bigcup(\mathbb{H}),$
a subspace $V\subseteq\mathbb{H}$ is invariant or a $G-invariant\
subspace$ if $\rho(g)V\subseteq V$ for all $g\in G.$ In this case
the restriction $\rho_V : G\rightarrow V$
  given by restricting $\rho(g)$ to $V$ is also a representation. If $\rho$
has no $G-invariant$ subspace other than $\{0\}$ and $\mathbb{H},$
then $\rho$ is irreducible. Equip $\mathbb{H}$ with an inner product
$\langle\cdot,\cdot\rangle,$ and define a new inner product
$\langle\cdot,\cdot\rangle^*$ which is preserved under the action
  of $\rho $ as follows:
  $$\langle v,w\rangle^*=\frac{1}{|G|} \sum\limits_{g\in G}\langle\rho(g)v,\rho(g)w\rangle,\
     v,w\in\mathbb{H}.$$
When $\rho$ is  not irreducible, there is a nontrivial invariant
subspace $V\subset \mathbb{H};$ and as $\langle\cdot,\cdot\rangle^*$
is invariant under each unitary map $\rho(g),$ the subspace
$V^\perp= \{u| v\in V, \langle u,v\rangle^*=0\}$ is also invariant.
Corresponding
 to the decomposition $\mathbb{H}=V\bigoplus V^\perp$, $\rho (g)$ has a natural decomposition
  $\rho(g)=\rho_V(g)\bigoplus \rho_{V^\bot}(g)$
for any $g\in G.$ Repeating this process, we see $\mathbb{H}$
  has a direct-sum decomposition: $\mathbb{H}=V_1\bigoplus\cdots\bigoplus V_k,$ and
  the following theorem holds.\\

\noindent{\bf Theorem (Complete reducibility) $2.1.$ } Any representation $\rho$ can be decomposed into
irreducible representations: $\rho=\sigma_1\oplus\cdots\oplus\sigma_k,$ where $\sigma_i=\rho_{V_i}.$\\

 \noindent{\bf Definition $2.1.$} Any two representations $\rho_1: G\rightarrow \bigcup(\mathbb{H}_1)$ and
 $\rho_2: G\rightarrow \bigcup(\mathbb{H}_2)$ are equivalent if there is an isomorphism
 $\theta : \mathbb{H}_1\rightarrow \mathbb{H}_2$ such that
 $$\theta(\rho_1(g)v)=\rho_2(g)\theta (v),\ g\in G,\ v\in \mathbb{H}_1.$$\\

 \noindent{\bf Theorem $2.2.$} Any finite group $G$ has only a finite number of irreducible representations
 up to equivalence.\\

 Let $G^*$ denote a set of representations containing exactly one from each equivalence class.
 There are two important representations in our analysis: one is the trivial representation
 $1: G\rightarrow \bigcup ( \mathbb{C} [G] ),\  i.e.,\ g\rightarrow id,$
 which is irreducible; the other is the regular representation
   $R : G\rightarrow \bigcup (\mathbb{C} [G] ),$ where $\mathbb{C}[G]$ is the vector space
 over the complex field $\mathbb{C}$ generated by $G,$ that is
 $$\mathbb{C} [G]=\left\{\left.\sum\limits_{g\in G}\alpha_g.g\right\vert \alpha_g \in \mathbb{C}\right\}.$$
Notice $(Rg)\left(\sum \alpha_h h\right)=\Sigma\alpha_hRgh;$ and $R$
is not irreducible and
$$R=\bigoplus\limits_{\rho\in
    G^*}\underbrace{\rho\oplus\cdots\oplus \rho}\limits_{d_\rho}.$$
Since $\mathbb{C}[G] $ has dimension $|G|,$ we have $|G|=\sum
d_\rho^2.$ Define $D=D(G)$ by $D=\sum d_\rho.$ Then
   $$\sqrt{|G|}<D(G)\leq |G|.$$

\section{Proof of Theorem 1.8}
Firstly, we verify the conditions in the proof of [6] Theorem 5 are
satisfied. Note [6] Theorem 5 is stated as Theorem 1.4 in our paper.

 Let $s_1,\cdots,s_k$ be chosen independently and uniformly at random from
 $G$ and
     $$S=\{s_1,\cdots,s_k\}.$$
     From theorem 1.6, we know that if  $\lambda_1>\lambda_2,$
then for any $0<\alpha\leq \frac{m}{n},$ $H$ is an
$\left(n,\frac{m}{n},l,\alpha,c(\alpha)\right)~bsc$ with
$$ c(\alpha)\geq\frac{l^2}{[\alpha(lr-\lambda_2)+\lambda_2]}.$$
Moreover the eigenvalue ordering of $X\left(G,H;SS^{-1}-\{ke\})\right)$ has same magnitude relation to the eigenvalue ordering of $X\left(G,H;SS^{-1} )\right).$ Thus we consider  $X\left(G,H;SS^{-1}-\{ke\})\right)$ following.

Let $B=SS^{-1}$ be the multiset as the product of $S$ and $S^{-1},$
and  $B^*=SS^{-1}-\{ke\},$ where $e$ is the unit element of $G,$
$$s=\frac{1}{2(k^2-k)}\sum\limits_{a\in B^*}(a+a^{-1}).$$
Notice the matrix of the linear operator
$$R(s)=\frac{1}{2(k^2-k)}\sum\limits_{a\in B^*}(R(a)+R(a^{-1}))$$
with respect to the standard basis of $\mathbb{C}[G]$ is the
normalized adjacency matrix of $X(G,SS^{-1}-\{ke\}),$ and its eigenvalue
$1$ corresponds to the trivial representation. By the decomposition
of $R,$ we have
 $\mu^*=\max \limits_{\rho \in G^*\backslash \{1\}}\|\rho(s)\|,$ where $\|\cdot\|$
is the operator norm.

For any  non-trivial representation $\rho$ of $G,$ let
  $$B_{a}=\frac{1}{2}\left[\rho (a)+\rho\left(a^{-1}\right)\right],\ a\in
     B^*,$$
and $\mu_1,\cdots,\mu_{d_\rho}$ be the eigenvalues of the
$\sum\limits_{a\in B^*} B_{a}$ arranged in decreasing order of their
absolute values, and $\lambda$ be an eigenvalue of
$\sum\limits_{a\in B^*}B_{a}$ chosen uniformly at random. Then
$$\mathbb{E}\left[\lambda \vert
    s_1,\cdots,s_k\right]=\frac{1}{d_\rho}\sum\limits_{a\in B^*}
   T_r(B_{a}),\ \mathbb{E}\left[\lambda\right]=\frac{1}{d_\rho}\mathbb{E}\left[\sum\limits_{a\in B^*}
   T_r(B_{a})\right].$$
By the decomposition of $R,$ $\sum\limits_{g\in G}\rho(g)$ is the
zero operator, and further we have
\begin{eqnarray*}
\mathbb{E}[\lambda]&=&\frac{1}{d_{\rho}}\sum\limits_{1\leq i\neq j\leq
   k}\mathbb{E}\left[Tr\left(B_{s_is_j^{-1}}\right)\right]
   =\frac{k^2-k}{d_{\rho}}\frac{1}{\vert G\vert^2 }\sum\limits_{x, y\in G}Tr\left(B_{xy^{-1}}\right)\\
&=&\frac{k^2-k}{2d_\rho\vert
   G\vert^2 }\sum\limits_{x, y\in
   G}Tr\left(\rho\left(xy^{-1}\right)+\rho\left(yx^{-1}\right)\right)\\
&=&\frac{k^2-k}{2d_\rho\vert
   G\vert^2 }\sum\limits_{y\in
   G}Tr\left(\sum\limits_{x\in
   G}\left(\rho\left(xy^{-1}\right)+\rho\left(yx^{-1}\right)\right)\right)\\
&=&\frac{k^2-k}{d_\rho\vert
   G\vert^2 }\sum\limits_{y\in
   G}Tr\left(\sum\limits_{g\in
   G}\rho\left(g\right)\right)=0.
\end{eqnarray*}
Now let $\lambda_i=\mathbb{E}\left(\lambda \vert B_{a_1},\cdots,B_{a_i} \right),$ then $\lambda_0,\cdots,\lambda_{k^2-k }$ is a martingale with $\vert \lambda_i-\lambda_{i-1}\vert \leq  1 ,$ since $\lambda_i-\lambda_{i-1}=\frac{Tr(B_{a_i})}{ d_\rho}$ and $\rho (a_i), \rho\left(a_i^{-1}\right)$ are unitary operators, thus $\vert Tr(B_{a_i})\vert\leq d_\rho.$ Applying Hoeffding-Azuma inequality, we conclude that $$Pr\left( \vert \lambda \vert\geq \varepsilon \left(k^2-k\right)\right)\leq 2exp\left\{-\left(k^2-k\right)H_{\frac{1}{2}}\left(\frac{1}{2}+\varepsilon\right)\right\}$$
At last we have that
\begin{eqnarray*}
 Pr \left( \left\|\rho\left(s\right)   \right \|\geq \varepsilon\right)
&=&Pr \left( \left\|\frac{1}{k^2-k} \left( \sum_{i=1}^{k^2}B_{a_i} \right) \right \|\geq \varepsilon\right)\\
&=&Pr \left(  \left |   \lambda_{1 }     \right |\geq  \varepsilon \left(k^2-k\right)\right)\\
 &\leq& \sum_{i=1}^{d_\rho}Pr\left( \left |   \lambda_{i }    \right |\geq  \varepsilon \left(k^2-k\right)\right)\\
&=& d_\rho Pr\left( \left| \lambda \right|\geq  \varepsilon \left(k^2-k\right)\right)\leq 2d_{\rho}exp\left\{-\left(k^2-k\right)H_{\frac{1}{2}}\left(\frac{1}{2}+\varepsilon\right)\right\}\\
\end{eqnarray*}Therefore, summing over all irreducible non-trivial representations of $G,$ we have the following: For any $0<\varepsilon<1,$
$$P_r\left(\mu^*(X(G,A^*)) > \varepsilon\right)\leq 2D\exp \left\{-(k^2-k)
   H_ \frac{1}{2}\left(\frac{1}{2}+\varepsilon\right)\right\}.$$\\

Notice $  AA^T - k I $ is the adjacency matrix of
$X\left(G,H; B^*\right),$ and the matrix
$$M_L^*=\frac{1}{2\left(|S|^2-k\right)}(AA^T -kI)$$ is the normalized adjacency matrix of
$X\left(G,H;B^*\right).$ If  $ \lambda_1^* \geq   \lambda_2^* \geq \cdots\geq  \lambda_n^* $ is the ordered normalized eigenvalues of $AA^T-kI,$ then  $\frac{\left(2\left(\vert S\vert^2-k\right)\right)\lambda_1^*+k}{2\vert S\vert^2}\geq \frac{\left(2\left(\vert S\vert^2-k\right)\right)\lambda_2^*+k}{2 \vert S\vert^2}\geq \cdots\geq \frac{\left(2\left(\vert S\vert^2-k\right)\right)\lambda_n^*+k}{2\vert S\vert^2}$ are the ordered normalized eigenvalues of $ AA^T.$ Thus we can obtain that if $\lambda_1^*> \lambda_2^*,$ then  $\frac{(2\vert S\vert^2-k)\lambda_1^*+k}{2\vert S\vert^2}>\frac{(2\vert S\vert^2-k)\lambda_2^*+k}{2\vert S\vert^2}.$ Let $M_L=\frac{1}{2|S|^2 }AA^T,$ then by Theorem 1.5, we see that for any $0<\varepsilon<1,$  $$P_r\left(\mu^* \left(M_L\right)>\varepsilon+\frac{1-\varepsilon}{2k}\right)\leq
   2D(G,H) exp\left\{-(k^2-k) H_{\frac{1}{2}}\left(\frac{1}{2}+\varepsilon \right)\right\}.$$\\

$\hfill \Box$\\

\noindent{\bf Corollary $3.1.$} If $\vert S\vert =k$ is large
enough, then the random bi-coset graph $X=C(G,L,N;S)$ of
$G$ with respect to $L,N$ and $S$ is almost always a concentrator.\\

\noindent{\bf Remark $3.1.$} Let $S^*=S\cup \{1\}.$ Then the random
bi-coset graph $X=C\left(G,L,N;S^{*}\right)$ of $G$ with respect to
$L,N$ and $S^*$ is almost always a concentrator if $\vert S\vert =k$
is large enough.

\section{Bsc and semi-vertex transitive graph from generalized polygons}
  The generalized polygons are incidence structures consisting of points and lines
of which the bipartite graphs have diameter $D$ and girth $2D$ for
some integer $D$. Tanner$\cite{18}$ proved every bipartite graph of
the generalized polygons whose any point is incident on $s+1$ lines
and any line is incident on $r+1$ points, is a good $bsc$ for
$D-gons$ with $D=3,4,6,8.$

    Note that for $D=6$, the smallest thick generalized hexagon has order
    (2,2), namely $s=2,r=2.$
Let $X$ be the bipartite graph of the generalized hexagon of order
(2,2). It is known that any automorphism of $X$ does not act
transitively on the vertices, and $X$ has two orbits that are the
two halves of the bipartition, i.e. $X$ is semi-vertex transitive.
For other $D$ and $(s,r),$ we do not know generally whether the
related bipartite graphs are semi-vertex transitive; but there are
some $D$ and $(s,r),$ the related bipartite graphs are not
semi-vertex transitive.\\

\noindent{\bf Proposition $4.1$.} $X$ is a semi-vertex transitive
$bsc.$

\section{Bsc and semi-vertex transitive graph from designs of the Mathieu groups}
  \noindent{\bf Definition $5.0.1$.}
Given any natural numbers $v,k,\gamma,t$ with $v>k>t\geq 2.$ Let $P$
be a set of elements called points and $B$ a set of subsets of $P$
called blocks such that
\begin{eqnarray*}
&&(a)\ |P|=v,\ (b)\ |B_i|=k~\mbox{for\ all}\ B_i\in B,\\
&&(c)\ \mbox{any}~ t-tuple~ \mbox{of \ points\ is\ contained\  in\
   exactly}\
      \gamma\ \mbox{blocks}.
\end{eqnarray*}
Then the system $\mathcal{D}=(P,B)$ is called a $t-(v,k,\gamma)$ design.\\

\noindent{\bf Proposition $5.0.1^{\cite 8}$. }Assume $\mathcal{D} = (P, B)\ \mbox{is\ a}\
 t-(v, k, \gamma)$ design. Then for any natural number $s\leq t$ and any
 subset $S$ of $P$ with $|S| = s,$ the total number of blocks incident with each element of $S$ is
given by
$$ \gamma_s=\gamma \frac{(v-s)(v-s-1)\cdots(v-t+1)}{(k-s)(k-s-1)\cdots(k-t+1)}.$$
Particularly, a $t-(v,k,\gamma)\ \mbox{design\ is\ also\ an}\ s-(v, k, \gamma_s )$ design.\\

 Let $r := \gamma_1$ be the total number of blocks incident with a given point.\\

\noindent{\bf Proposition $5.0.2^{\cite 8}$. }Let $
\mathcal{D}=(P,B)$ be a $t-(v,k,\gamma)$ design. Then
$$bk=vr,\  {v\choose t} \gamma=b {k\choose t},\ r(k-1) = \gamma_2(v-1).$$

\noindent{\bf Definition $5.0.2$.} A balanced incomplete block
design ($BIBD$) with parameters $(v,b,r,k,\lambda)$ is an arrangement
of $v$ distinct objects $b$ blocks such that each block contains
exactly $k$ distinct objects, each object occurs in exactly $r$
different blocks, and every pair of distinct objects $a_i,\ a_j$
occur together in exactly $\lambda$ blocks. Obviously a $BIBD$ is a
$2-(v,b,r,k,\gamma_2)$ design.\\

\noindent{\bf Definition $5.0.3$.} Assume $P=\{a_1,\cdots,a_v\}$ is
a set of $v$ objects, and $B=\{B_1,\cdots,B_b\}$ is a set of $b$
blocks consisting of elements of $P,$ and $\mathcal{D}=(P,B)$ is a
$BIBD$ with parameters $(v,b,r,k,\lambda).$ Define $X_{\mathcal{D}}$
as the following bipartite graph: This bipartite graph has $A$ as
the set of left vertices and $B$ as the set of right vertices such
that there is an edge between
 an $a_i$ and a $B_j$ with $a_i\in B_j.$

Let $v\times b$ matrix $A=(a_{i,j})$ be the incidence matrix of
$\mathcal{D},$ where $a_{i,j}=I_{\{a_i\in B_j\}}.$\\

%
%
%

 \noindent{\bf Lemma $5.0.1^{\cite{7}}$.} The character and minimal polynomials of $AA^T$ are
respectively
\begin{eqnarray*}
 &&P(x)=(x-rk)(x-(r-\lambda))^{v-1},\ x\in\mathbb{R}^1,\\
 &&m(x)=(x-rk)(x-(r-\lambda)),\ x\in\mathbb{R}^1.
\end{eqnarray*}
Hence $\lambda_1=rk,\ \lambda_2=r-\lambda~\mbox{and} ~r-\lambda<rk,$
and $X_{\mathcal{D}}$ is a $bsc.$\\

  For any graph $X$ with $n$ vertices, its adjacency matrix has $n$
eigenvalues which are denoted by
$\mu_0\geq\mu_1\geq\cdots\geq\mu_{n-1}$ with
decreasing order.\\

\noindent{\bf Definition $5.0.4$ (Bipartite Ramanujan Graphs).} For
any $(c,d)-\mbox{regular\ bipartite\ graph}$ $X,$ we call it a Ramanujan graph if $\mu_1\leq \sqrt{c-1}+\sqrt{d-1}.$\\

Recall H$\phi$holdt and Janwal$\cite{7}$ proved the bipartite graph
of a $(v,b,r,k,\lambda)$ $BIBD$ is a $(r,k)-regular$ bipartite
Ramanujan graph with $\mu_1=\sqrt{r-\lambda},$ and it is the optimal
expander graph with these parameters.\\

 Let $\mathcal{D}=(P,B)$ be a $t-(v,k,\gamma)$ design and $p$ a point of $\mathcal{D}.$
Define a new design $\mathcal{D}_p$ depending on $p$ as follows: its
point set is $P-\{p\}$, and its block set is
$$\{B-\{p\}\ |\ B\ \mbox{is\ any\ block\ containing}\ p\ \mbox{in} ~\mathcal{D}\}.$$
Note $\mathcal{D}_p$ is a $(t-1)$-$(v-1,k-1,\gamma)$ design and we
call it the contraction of ${\mathcal{D}}$ at $p.$

Given a permutation $g$ of $P,$ if for any block $B$ of
$\mathcal{D},$ $B^g=\{x^g:=g(x)\ |\ x\in B\}$ is also a block of
$\mathcal{D},$ then $g$ is called a automorphism of $\mathcal{D}.$
Clearly all automorphisms of $\mathcal{D}$ forming a group, which is
called the automorphism group of $\mathcal{D}$ and denoted by
$Aut(\mathcal{D}).$\\

In Subsections 5.1-5.2, we will construct some semi-vertex
transitive $bscs$ from Mathieu groups.\\

\subsection{The large Mathieu groups}
  The projective plane $PG(2,4)$ can be extended 3 times leading to the unique designs with
parameters $3-(22,6,1),$ $4-(23,7,1)$ and $5-(24,8,1),$ denoted by
$\mathcal{M}_{22},\mathcal{M}_{23}$ and $\mathcal{M}_{24}$
respectively. It is well-known that $Aut(\mathcal{M}_{i})$ is the
$Mathieu ~group~M_i$ for any $i\in \{22,23,24\}.$

Let $\Omega $ be a set with $|\Omega|=n,$ and $H$ a subgroup of the
symmetric group $Sym(n).$ We say $H$ is transitive on $\Omega$ if
for every $a$ and $b$ in $\Omega,$ there exists $\pi\in H$ such that
$a^\pi=b.$ Given a natural number $k\leq n.$ If for every list of
$k$ distinct points $a_1,\cdots,a_k$ and very list of $k$ distinct
points $b_1\cdots,b_k,$ there exists a $\pi\in H$ such that
$a_i^\pi=b_i$ for all $i,$ then $H$ is said to be $k-transitive.$
Clearly, if $H\neq \{1\}$ is $k-transitive,$ then it is
$m-transitive$ for all $m\leq k.$\\

\noindent{\bf Proposition $5.1.1^{\cite{14}}$.} The group $M_{24}$
is 5-transitive on the 24 points of $\mathcal{M}_{24},$ and the
group $M_{24-i}$ is $(5-i)$ transitive on $24-i$ points of
$\mathcal{M}_{24-i}$ for any $i\in\{1,2\}.$ Particularly, $M_j$ is
transitive on the $j$ points for any $j\in\{24,23,22\}.$\\

\noindent{\bf Theorem $5.1.1$.} For any $i\in\{22,23,24\},$ let
$X_{\mathcal{M}_i}$ be the corresponding bipartite graph of the
design $\mathcal{M}_i.$ Then each $X_{\mathcal{M}_i}$ is a
semi-vertex
transitive $bsc.$\\

\noindent{\bf Proof:} Fix an $i\in\{22,23,24\}.$ Let
$\mathcal{M}_i=(P,B).$ Then for any $a ~\mbox{and} ~b\in P,$ there
exists $g\in M_i$ with $a^g=b.$ Assume $a\in B_j$ and $B_j\in B,$
since $g$ is a automorphism of the design, we have that $B_j^g\in
B,$ and $g$ maps the neighbors of $a$ to the neighbors of $b,$ and
$X_{\mathcal{M}_i}$ is semi-vertex transitive.

By Proposition 5.0.2, $3-(22,6,1),~4-(23,7,1)$ and $5-(24,8,1)$ are
$(22,77,21,6,5),$ $(23,506,77,6,21)$ and $(24,759,253,8,77)$ $BIBDs$
respectively. From Lemma 5.0.1 and Theorem 5.1.1, we obtain that
each $X_{\mathcal{M}_i}$ is a $bsc.$  Moreover, each
$X_{\mathcal{M}_i}$ is a bipartite Ramanujan graph by the main
result of H$\phi$holdt and Janwal$\cite{7}.$ $\hfill \Box$

\subsection{The small Mathieu groups}

Let $\Omega=\{1,2,\cdots,12\}$ and consider the following
permutations on $\Omega:$
\begin{eqnarray*}
\mu &=&(1 ~2~ 3)(4~ 5~ 6)(7~ 8 ~9),\ \ \ \ \ \ a=(2 ~4 ~3~ 7)(5~ 6~ 8~ 9),\\
b&=&(2 ~5~ 3~ 9)(4~ 8~ 7~ 6),\ \ \ \ \ \ \ \ \ \ x=(1 ~10)(4~ 5)(6~ 8)(7~ 9),\\
y&=&(10~ 11)(4~ 7)(5~ 8)(6~ 9),\ \ \ \ z=(11~ 12)(4~ 9)(5 ~7)(6~ 8).
\end{eqnarray*}
Then $M_{12}=\langle\mu,a,b,x,y,z\rangle.$ Let
$\Delta=\{1,2,3,10,11,12\}$ and $\mathcal{D}_{12}=(\Omega,B),$ where
$$B=\{\Delta^g|g\in M_{12}\}.$$
Note $M_{12}$ is the automorphism group of the design
$\mathcal{D}_{12},$  which is a $5-(12,6,1)$ design. Let $p$ be any
element of  $\Omega,$ say $p=12.$ Then the contraction
$(\mathcal{D}_{12})_p$ of ${\mathcal{D}_{12}}$ at $p=12$ is a
$4-(11,5,1)$ design and
$M_{11}=Aut\left(\left(\mathcal{D}_{12}\right)_p\right).$ Similarly
we can construct a $3-(10,4,1)$ design and a $2-(9,3,1)$ design,
$M_{10}$ and $M_9$ are their automorphism groups respectively. So we
can denote the designs $4-(11,5,1),3-(10,4,1),2-(9,3,1)$ by
$\mathcal{D}_{11},\mathcal{D}_{10},\mathcal{D}_{9}$ respectively.

Notice $5-(12,6,1),4-(11,5,1),3-(10,4,1),2-(9,3,1)$ are respectively
$(12,132,66,6,30),$ $(11,66,30,5,12),(10,30,12,4,4),(9,36,8,3,1)$
$BIBDs$. Similarly to Theorem 5.1.1, we can prove\\

\noindent{\bf Theorem $5.2.1$.} Let $X_{\mathcal{D}_i}$ be the
corresponding bipartite graph of the design $\mathcal{D}_i,$ then
for any $i\in\{9,10,11,12\},$ $X_{\mathcal{D}_i}$ is a semi-vertex
transitive $bsc.$\\

 \section{Symmetric  group and a sequence of concentrators}
Assume $X=(V,E)$ be a connected graph with $n$ vertices and $Q=Q_X=diag(d(v))-A_X$ where $A_X$ is its adjacency matrix. Let $\lambda(X)$ be the second smallest eigenvalue of $Q.$ When $X$ be a $d-regular$ graph, then $\lambda(X)$ is the difference between $d$ and the second largest eigenvalue of $X.$ The following result holds.
\\

\noindent{\bf Theorem $6.1^{\cite{1}}$.} If a $d$-regular graph
$X=(V,E)$ is an $(n,d,\varepsilon)-magnifier,$ then
    $$\lambda(G)\geq \frac{\varepsilon^2}{4+2 \varepsilon^2}.$$\\

Recall the following result from Kassabov$\cite{9}.$\\

\noindent{\bf Theorem $6.2^{\cite{9}}$.} For every natural number
$n,$ there is a generating set $S_n$ (of size at most $\ell$) of the
alternating group $Alt(n)$ such that the Cayley graphs
$X\left(Alt(n), S_n\right)$ form a family of $\varepsilon-expanders\
(magnifiers).$ Here $\ell$ and $\varepsilon > 0 $ are some universal
constants. Similarly there is a generating set $\widetilde{S}_n$ of
the symmetric group $Sym(n)$ with the same property.\\

 Let $L$ be a subgroup of $Alt(n)$ or $Sym(n).$ Then the bi-coset graphs
 $$X=C\left(G,L,\{1\},\left(S_n\cup\{1\}\right)L\right)\ \mbox{and}\
    Y=C\left(G,L,\{1\},\left(\widetilde{S}_n\cup\{1\}\right)L\right)$$
are semi-vertex transitive graphs, and are respectively
$\left(\left|\left(S_n\cup\{1\}\right)L\right|,\left|S_n\cup\{1\}\right|\right)-$bipartite
graph and
$\left(\left|\left(\widetilde{S}_n\cup\{1\}\right)L\right|,
   \left|\widetilde{S}_n\cup\{1\}\right|\right)-\mbox{bipartite\
   graph}.$\\

\noindent{\bf Theorem $6.3$.} The bi-coset graphs
$$X=C\left(G,L,\{1\},\left(S_n\cup\{1\}\right)L\right)\ \mbox{and}\
    Y=C\left(G,L,\{1\},\left(\widetilde{S}_n\cup\{1\}\right)L\right)$$
 are $bscs.$\\

\noindent{\bf Proof:} We only prove the case for alternating groups.
Let $M_L^*=AA^T,$ where $A$ is the adjacency matrix of $X.$ Then
$M_L^*$ is a symmetric nonnegative definite matrix and has
nonnegative eigenvalues, which can be considered as the adjacency
matrix of the coset graph
$$D=\left(G/L,(S_n\cup\{1\})(S_n\cup\{1\})^{-1}\right).$$
Notice $D$ is a
$2\left|(S_n\cup\{1\})(S_n\cup\{1\})^{-1}\right|-regular$ graph.
Assume
 $$G/L=\{L_{g_1},\cdots,L_{g_m}\}\ \mbox{and}\
   A=\{L_{g_{i1}},\cdots,L_{g_{ir}}\}\ \mbox{with}\ r\leq\left\lfloor\frac{m}{2}\right\rfloor.$$
Then by $X(Alt(n), S_n)$ is an $\varepsilon-expander\ (magnifier),$
we have
$$\left|A\left((S_n\cup\{1\})(S_n\cup\{1\})^{-1}\right)\backslash
A\right|\geq |A(S_n\cup\{1\})\backslash A|=|AS_n\backslash A|\geq
\varepsilon |A|.$$ Hence we obtain
\begin{eqnarray*}
&&\left|\left(A\left((S_n\cup\{1\})(S_n\cup\{1\})^{-1}\right)/L\right)\backslash
  (A/L)\right|\\
&&\ \ \ \  \geq \left|\left(A(S_n\cup\{1\})/L\right)\backslash
      (A/L)\right|=\left|\left((AS_n)/L\right)\backslash
    (A/L)\right|\geq \varepsilon \left|A/L\right|.
\end{eqnarray*}
Therefore,
$$D=\left(G/L,(S_n\cup\{1\})(S_n\cup\{1\})^{-1}\right)$$
is an $\varepsilon-expander$ graph, namely an
$(n,d,\varepsilon)-magnifier.$ By Theorem 6.1,
$$X=C\left(G,L,\{1\},(S_n\cup\{1\})L\right)$$ is
a $bsc.$ $\hfill\Box$

\section{Concluding remarks}
We prove the random bi-coset graphs are almost always concentrators,
and construct some examples of semi-vertex transitive $bscs.$
Because generalized $D-gons$ do not exist for arbitrary parameters
$n$ and $k$, Tanner$\cite{18}$ did not provide a complete solution
to the problem of constructing concentrators. However, we can get a
sequence of concentrators by symmetric groups or alternating groups
with their appropriate subgroups for arbitrary parameter $n$ and
some $k\leq \ell$ in Section 6.

\end{document}